\documentclass{amsart}
\usepackage{amsfonts}
\usepackage{amsmath}
\usepackage{amssymb}
\usepackage[sans]{dsfont}

\usepackage{graphicx}
\usepackage{float}

\newtheorem{theorem}{Theorem}[section]
\newtheorem{lemma}[theorem]{Lemma}
\theoremstyle{definition}
\newtheorem{definition}[theorem]{Definition}
\newtheorem{example}[theorem]{Example}

\theoremstyle{proposition}
\newtheorem{proposition}[theorem]{Proposition}
\theoremstyle{remark}
\newtheorem{remark}[theorem]{Remark}
\theoremstyle{corollary}
\newtheorem{corollary}[theorem]{Corollary}

\numberwithin{equation}{section}

\graphicspath{%
    {converted_graphics/}
    {/}
}
\begin{document}
\def\R{{\mathbb R}}
\def\C{{\mathbb C}}
\def\N{{\mathbb N}}
\def\DD{{\mathbb D}}
\def\diam{\mbox{\rm diam}}
\def\rr{{\cal R}}
\def\e{\emptyset}
\def\dQ{\partial Q}
\def\dk{\partial K}
\def\endofproof{{\rule{6pt}{6pt}}}
\def\di{\displaystyle}
\def\dist{\mbox{\rm dist}}
\def\u-{\overline{u}}
\def\du{\frac{\partial}{\partial u}}
\def\dv{\frac{\partial}{\partial v}}
\def\dt{\frac{d}{d t}}
\def\dx{\frac{\partial}{\partial x}}
\def\con{\mbox{\rm const }}
\def\Box{\spadesuit}
\def\ii{{\bf i}}
\def\curl{{\rm curl}\,}
\def\dive{{\rm div}\,}
\def\grad{{\rm grad}\,}
\def\dist{\mbox{\rm dist}}
\def\pr{\mbox{\rm pr}}
\def\pp{{\cal P}}
\def\supp{\mbox{\rm supp}}
\def\Arg{\mbox{\rm Arg}}
\def\In{\mbox{\rm Int}}
\def\Re{{\rm Re}\:}
\def\Im{{\rm Im}\:}
\def\li{\mbox{\rm li}} 
\def\ep{\epsilon}
\def\tr{\tilde{R}}
\def\be{\begin{equation}}
\def\ee{\end{equation}}
\def\cn{{\mathcal N}}
\def\sn{{\mathbb  S}^{n-1}}
\def\Ker {{\rm Ker}\:}
\def\el{E_{\lambda}}
\def\Rc{{\mathcal R}}
\def\Ha{H_0^{ac}}
\def\la{\langle}
\def\ra{\rangle}
\def\Ko{\Ker G_0}
\def\Kd{\Ker G_b}
\def\hc{{\mathcal H}}
\def\caH{\hc}
\def\caO{{\mathcal O}}
\def\Hp{\hc_b^{\perp}}
\def\Sc{{\mathcal S}}
\def\Kc{{\mathcal K}}
\def\Dp{D_{+}^{\rho}}
\def\Dm{D_{-}^{\rho}}
\def\D2p{D_{+}^{\2\rho}}
\def\D2m{D_{-}^{\2\rho}}
\def\12{\frac{1}{2}}
\def\S{{\mathbb S}}
\def\one{\mathds{1}}

\title[Spectral problems] {Spectral problems for non elliptic symmetric systems with dissipative boundary conditions}


\author[F. Colombini, V. Petkov] {Ferruccio Colombini, \ Vesselin Petkov}
\address{Dipartimento di Matematica, Universit\`a di Pisa, Italia}
\curraddr{}
\email{colombini@dm.unipi.it}
\address{Institut de Math\'ematiques de Bordeaux, 351,
Cours de la Lib\'eration, 33405  Talence, France}
\email{petkov@math.u-bordeaux1.fr}

\thanks{}

\author[J. Rauch]{\ Jeffrey Rauch}
\address{Department of Mathematics, University of Michigan, USA}
\email{rauch@umich.edu}
\curraddr{}
\thanks{JR was partially supported by NSF grant
DMS 0807600
}
\thanks{VP was partially supported by ANR project Nosevol BS01019 01}

\date{}

\dedicatory{}


\begin{abstract} 
This paper considers and extends  spectral and scattering theory 
to  dissipative 
symmetric systems that may have zero speeds and  in particular
to strictly dissipative boundary conditions for Maxwell's equations. 
Consider symmetric systems $\partial_t - \sum_{j=1}^n A_j \partial_{x_j}$ in $\R^n,\: n \geq 3$, $n$ odd,
in a smooth connected exterior domain 
$\Omega :=\R^n \setminus \bar{K}$.
Assume that  the rank of $A(\xi) = \sum_{j= 1}^n A_j \xi_j$ is constant for $\xi \not= 0.$ For 
maximally dissipative boundary conditions on $\Omega :=\R^n \setminus \bar{K}$ with bounded open domain $K$ the solution of the boundary problem in $\R^{+} \times \Omega$ is described by a contraction semigroup $V(t) = e^{t G_b},\:t \geq 0.$ Assuming 
coercive  conditions for $G_b$ and its adjoint $G_b^*$
on the complement of their kernels, we prove that  the spectrum of $G_b$ in the open half plane $\Re z < 0$
 is formed only by isolated eigenvalues with finite multiplicities. 
 \end{abstract}

\maketitle

 

\section{Introduction}

This paper is devoted to dissipative symmetric hyperbolic
systems on exterior domains.  The systems studied
have two properties that render their study difficult.
First, they are dissipative and not conservative so self
adjointness techniques are not available.  Second, the symbol $A(\xi)$ is not
elliptic. Third, we may have characteristic roots with variable multiplicities.  However, we suppose that
the failure of ellipticity is uniform in $\xi\ne 0$ and that
coercivity estimates hold on the orthogonal complement of the space $H_b$ spanned by
 the stationary solutions.  
We prove that
the spectrum of the generator $G_b$  restricted to $H_b^{\perp}$ 
 consists of the imaginary axis plus at most discrete subset of ${\rm Re}\,z<0$. Moreover, for the Maxwell system we show that 
 the discrete part of the spectrum is stable under perturbations.
For the wave equation with dissipative boundary conditions  such results have been proved 
before (see \cite{LP2}).

Suppose that $n\ge 3$ is odd, and that 
 $K \subset \R^n$ is a bounded open set with smooth boundary. Assume that the exterior domain $\Omega :=\R^n \setminus \bar{K}$ is connected.
 Choose $\rho>0$ so that $K\subset\{x \in \R^n:|x|\le \rho\}$.
 Suppose that for $j =1,...,n,$ $A_j$ are $(r\times r)$ symmetric  matrices and set $A(\xi) = \sum_{j=1}^n A_j \xi_j$. Assume that 
 Rank $A(\xi) = r - d_0 > 0$ is independent of $\xi \not= 0$. Define  $G := \sum_{j=1}^n A_j \partial_{x_j}$ and for $x\in \partial \Omega$, 
denote by $\nu(x) = (\nu_1(x),...,\nu_n(x))$ the unit outward normal to $ \Omega$,
and $A( \nu(x)) = \sum_{j=1}^n A_j \nu_j(x)$. 

Suppose that $\cn(x) \subset \C^r$ is a linear space depending smoothly on $x \in \partial \Omega$ such that

$(i)\:\:\:   \langle A(\nu(x)) u(x), u(x) \rangle \leq 0$ for all $u(x) \in \cn(x),$

$(ii) \:\:\:  \cn(x)$ is maximal with respect to (i).

Then the initial boundary value problem
\begin{equation} \label{eq:1.1}
\begin{cases} (\partial_t - G) u = 0 \quad{\rm in}\quad \R^+ \times \Omega,\\
u(t,x) \in \cn(x) \quad {\rm for}\quad t\geq 0,\quad 
x \in \partial \Omega,\\
u(0, x) = f(x) \quad
 {\rm in}\quad
  \Omega \end{cases}
\end{equation}
yields a contraction semigroup $V(t) = e^{tG_b},\: t \geq 0,$ in $\hc = L^2(\Omega: \C^r)$ with generator $G_b$. The domain 
$D(G_b)$ of the generator
 is the  closure with respect to the graph norm $(\|g\|^2 + \|Gg\|^2)^{1/2} $ of functions $g(x) \in C_{(0)}^1(\bar{\Omega}: \C^r)$ satisfying the boundary condition $g(x)\vert_{\partial \Omega} \in \cn(x)$.
 
Introduce the unitary group $U_0(t) = e^{t G_0}$ on $H_0 = L^2(\R^n: \C^r)$ 
solving  the Cauchy problem
\begin{equation} \label{eq:1.2}
\begin{cases} (\partial_t - G) u = 0 \quad{\rm in}\quad \R \times \R^n,\\
u(0, x) = f(x) \quad {\rm in}\quad \R^n.
 \end{cases}
\end{equation}

Define $\hc_b \subset \hc$ to be the space generated by the eigenvectors of $G_b$ with eigenvalues $\mu \in \ii\R$ and let $\Hp$ be the orthogonal complement of $\hc_b$ in $\hc$. The generator $G_0 =  \sum_{j=1}^n A_j \partial_{x_j}$ is skew self-adjoint in $H_0$ and the spectrum of  $G_0$ on the space $\Ha = ({\Ker}\: G_0)^{\perp} \subset H_0$ is absolutely continuous (see Chapter IV in \cite{P2}).

Throughout this paper we suppose that $G_b$ and $G_b^*$ satisfy the following coercive estimates

$(H):\:$ For each $f \in D(G_b) \cap (\Kd)^{\perp}$ we have 
\begin{equation} \label{eq:1}
\sum_{j=1}^n \|\partial_{x_j} f\| \leq C (\|f\| + \|G_b f\|)
\end{equation}
with a constant $C > 0$ independent of $f.$ \\

$(H)^*:\:$ For each $f \in D(G_b^*) \cap (\Ker G_b^*)^{\perp}$ we have 
\begin{equation} \label{eq:2}
\sum_{j=1}^n \|\partial_{x_j} f\| \leq C (\|f\| + \|G_b^*f\|)
\end{equation}
with a constant $C > 0$ independent of $f.$\\


Notice that  when
$G$ is elliptic, that is $d_0 = 0$,  these estimates hold exactly when the 
associated boundary value
problems satisfy Lopatinski conditions.   The interesting cases are when $G$ is  not elliptic and the estimates only hold
on the orthogonal 
complement
to the kernels.
Majda \cite{M} proved
$(H)$ and $(H^*)$ when
the constant rank hypothesis holds and in addition the following condition
is satisfied.

\vskip.2cm
$(E):$ There exists a first order $(l \times r)$ matrix operator $Q = \sum_{j = 1}^n Q_j \partial_{x_j}$ so that 
$$Q(\xi) A(\xi) = 0,\:\: \Ker Q(\xi) = {\rm Range}\: A(\xi),$$
where $Q(\xi) = \sum_{j= 1}^n Q_j \xi_j.$
\vskip.2cm

This says that one has an exact sequence
$$
\C^r\stackrel{A(\xi)} {\longrightarrow}\
\C^r \stackrel{Q(\xi)}
{\longrightarrow}\
\C^r\,.
$$
When the kernel of $A(\xi)$ has constant rank it was is mentioned in \cite{M},  that we can always 
choose a $(r \times r)$ matrix valued polynomial $Q(\xi) = \sum_{|\alpha| \leq m} A_{\alpha} \xi^{\alpha}$ with matrix coefficients $A_{\alpha}$ such that
$$(E_g): \Ker Q(\xi) = {\rm Range}\: A(\xi),\: \forall \xi \neq 0.$$
 Moreover,  it was remarked in \cite{M} that by using $(E_g)$, the results concerning coercivity have analogue by applying the techniques for systems elliptic in Agmon-Douglis-Nirenberg sense. However, in \cite{M} we can not find proofs of these two statements. Since the existence and the form of the symbol $Q(\xi)$ for which $(E_g)$ holds are important for our argument, for the sake of completeness we discuss the construction of $Q(\xi)$  in Section 4. 

On the other hand, if  $(E)$ holds, the symbol  $Q(\xi)$ is a linear
function of $\xi$ and the second order matrix operator $(-G^2 + Q^*Q)$ is strongly
 elliptic.

 \begin{example}
 Our principal motivation is the Maxwell equations 
 $$
 \partial_t E - \curl B = 0, \qquad \partial_t B + \curl E = 0,\qquad \dive E = 0,\qquad \dive B = 0.
 $$
For this system one can take
$$
Q\Bigl(\begin{matrix} E \\ B \end{matrix}\Bigr)
\ : =\
 \Bigl(\begin{matrix} \dive E\\ \dive B \end{matrix} \Bigr).
 $$
 Majda \cite{M} proves that the  Maxwell system with strictly dissipative boundary conditions 
 satisfies hypothesis $(H)$ and $(H)^*$.
 \end{example}
For dissipative symmetric systems some solutions can have  global energy decreasing exponentially as $t \to \infty$ and it is possible also to have disappearing solutions. The precise definitions are the following.

\begin{definition} 
The solution $u := V(t)f$ is  a disappearing solution (DS), if there exists $T > 0$ such that $V(t)f = 0$ for $t \geq T.$
\end{definition} 

\begin{definition} 
The solution  $u := V(t)f$ is an
asymptotically disappearing solution (ADS), if there exists $\lambda \in \C$ with $\Re \lambda < 0$ and $f \not= 0$ such that $V(t)f = e^{\lambda t}f$.of the motivation 
\end{definition}

An (ADS) is generated by  an eigenfunction  
$G_b f = \lambda f$ with $f \neq 0$ and $\Re \lambda < 0.$ The existence of at {\it least one} (ADS) implies the non completeness of the wave operators $W_{\pm},$ that is ${\rm Ran}\: W_{-} \neq {\rm Ran}\: W_{+}$ (see \cite{P3} and Section 4). 
The disappearing solutions perturb considerably the inverse back-scattering problems related to the leading singularity of the scattering kernel (for more details see \cite{P3}). Disappearing solutions for the wave equation have been studied by Majda  \cite{M1} and for the symmetric systems by Georgiev \cite{G2}. On the other hand, (ADS) for Maxwell's
equations with dissipative boundary conditions on the sphere $|x| = 1$ 
are constructed in \cite{CPR}. One of the motivations for the present work is to show that the latter
(ADS) persist under small perturbations.  To show that we must
examine the spectrum of $G_b$ in the half plane $\Re z < 0.$ 
The main result is the following

\begin{theorem} 
\label{thm:main}
Assume the hypothesis $(H),\: (H)^*$ fulfilled. Then the spectrum of $G_b$ restricted to  $\Hp$ is formed  by a discrete set in $\{z \in \C:\Re z < 0\}$ of eigenvalues with finite multiplicities and a continuous spectrum on $\ii \R.$ 
\end{theorem}

A result similar to Theorem 1.4 has been proved in \cite{LP2} for the semigroup $V_{w}(t) = e^{t G_w},\: t \geq 0,$ related to the mixed problem
\begin{equation} \label{eq:11}
\begin{cases} (\partial_t^2 - \Delta) w = 0\: {\rm in} \: \R^+ \times \Omega,\\
\partial_{\nu} w - \alpha(x) \partial_t w = 0\:\: {\rm on}\: \R^+ \times \partial \Omega,\\
(w(0, x), w_t(0, x)) = (f_1, f_2),\end{cases}
\end{equation}
where $\alpha \geq 0$ is a smooth function on $\partial \Omega.$ Let $V_0(t) = e^{t G_f}$ be the unitary group related to the Cauchy problem for the wave equation in $\R \times \R^n.$ The strategy of Lax-Phillips in \cite{LP2} was to prove that for $\Re z > 0$ the difference
$$(G_w - z)^{-1}  - (G_f - z)^{-1}$$
is a compact operator. 

For symmetric systems whose generators are not elliptic and for which coercivity holds only for data satisfying additional constraints, like the divergence free constraints for the  Maxwell equations, the coercivity estimates for  the generators 
hold on different spaces.  This prevents a direct application of the strategy in \cite{LP2}.  
We show that for $\Re z < 0$ the operator $G_b - z$ on  $\Hp$ is Fredholm. The  main difficulty is to prove that there exists at least one point $z_0, \: \Re z_0 <0,$ such that the index of  $G_b - z_0$ is equal to  zero. Thus we must exclude the case when the entire half plane
$\{ \Re z < 0\}$ is included in the spectrum of $G_b$ as in the case of a semigroup related to the one side shift (see Example 1 in \cite{LP2}).

Our proof uses scattering theory and the statement of Lax and Phillips \cite{LP2} that $z_0, \: \Re z_0 < 0,$ can be an eigenvalue of $G_b$ only if $w_0 = - \ii z_0$ is such that the scattering matrix $\Sc(w_0)$ at $w_0$ has non trivial kernel or $w_0$ is a pole of $\Sc(w_0)$ (see Section 3 for the notation).  However, this statement was established in \cite{LP2} under the assumptions that  neither  $G_b$ nor $G_b^*$ have eigenfunctions with eigenvalues in $\Re z < 0$ which are both incoming and outgoing (see Section 4). This property plays important  role in the scattering theory related to unitary groups. In the case of dissipative boundary problems we deal with the proof of this property is based on Theorem 4.2 which has independent interest. As a first corollary we obtain that $G_b$ has no outgoing eigenfunctions, while $G_b^*$ has no incoming eigenfunction with eigenvalues in $\Re z < 0$ (see Corollary 4.5) which is stronger that the assumption needed to apply the Lax-Phillips result mention above. A second corollary is the fact that if we have at least one eigenfunction of $G_b$ with eigenvalue in $\Re z < 0$, the wave operators $W_{\pm}$ are not complete, that is ${\rm Ran}\: W_{+} \neq {\rm Ran}\: W_{-}$ (see Corollary 4.6). The proof of Theorem 4.2 is technically difficult since we must work with boundary value problems with characteristic boundary and the Holmgren type argument can be applied only for systems with non-characteristic boundary. To pass from a first order system to a high order elliptic system we exploit the condition $(E_g)$.

The link between the scattering matrix and the spectrum of $G_b$ makes possible to use the property of $\Sc(z)$ in a neighborhood of 0 in order to show that there exists a small neighborhood ${\mathcal U}$ of 0 in $\C$ such that ${\mathcal U} \cap \{\Re z < 0\}$ is in the resolvent set of $G_b$. This idea has been applied in \cite{LP2}. For this purpose it is necessary to examine the scattering matrix in a neighborhood of 0. The difference with \cite{LP2} is that in the case with non elliptic systems with characteristic of variable multiplicities the representation of the scattering matrix is more complicated than that for the wave equation and we are going to apply the representation of ${\mathcal S}(z)$ established for non elliptic systems with dissipative boundary conditions in \cite{P3}. Another consequence of our approach is Proposition 4.11 saying that the eigenvalues of $G_b$ in $\Re z < 0$ cannot have a finite accumulation point on $\ii \R$. 

Finally, in  our case the space $H_b$ is {\bf always} infinite dimensional. This fact has been proved by Majda \cite{M} if we have the condition (E). In the general case we use the property $(E_g)$ and the proof is a trivial repetition of that in \cite{M}. On the other hand, it was proved in \cite{GS}, that every eigenvalues $\ii z$ of $G_b$ with $z \in \R \setminus \{0\}$ has a finite multiplicity and the eigenvalues $\ii z$ of $G_b$ with $z \in \R\setminus \{0\}$ could have accumulation points only at 0 and $\pm \ii \infty.$ 

The paper is organized as follows. In Section 2 we study the Fredholm operator $G_b - z,\: \Re z < 0,$ and the proof of Theorem 1.1 is reduced to showing that the index of  $G_b - z$ is equal to zero. 
 In Section 3 we collect some facts about the scattering operator $S$ and the scattering matrix ${\mathcal S}(z).$ By applying a representation of the scattering kernel, we prove that ${\mathcal S}(z)$ is invertible in a small neighborhood of 0. In Section 4 Theorem 4.2 and Corollaries 4.5 and 4.6. We complete the proof of Theorem 1.1 by using the invertibility of ${\mathcal S}(0).$ In Section 5 we establish a result showing that the (ADS) for the Maxwell system are stable under the perturbations of the operator, the boundary and the boundary conditions. Combining Theorem 5.1 with the construction of (ADS) for Maxwell system in \cite{CPR}, we conclude that (ADS) exist for domains which are not spherically symmetric.\\

\noindent
{\bf Acknowledgment.}
We thank  Jean-Fran\c{c}ois Bony who  showed  us the construction of the symbol $Q(\xi)$ in Section 4. We would also thank the referee for many remarks and suggestions concerning the previous version of the paper.

\section{The Fredholm operator $G_b - z$}

Define $\hc_0 := (\Ker \,G_b)^{\perp}.$ For $\Re\: z < 0$, it is clear that $G_b - z:\: \hc_0 \cap D(G_b) \longrightarrow \hc_0.$

\begin{proposition} For every $R>\rho$ there is a constant $C_R > 0$ so that
for  every $u \in \hc_0 \cap D(G_b)$ and $\Re\:z < 0$,
\begin{equation} \label{eq:4}
\|u\|_{H^1(\Omega)} \leq C_R \Bigl(|z|(1  + \frac{1}{|\Re z|})\Bigr) ( \|(G_b - z)u\| + \|\one_{|x| \leq 2R} u \|).
\end{equation}

\end{proposition}
{\it Proof.}  The coercive estimate yields
$$
\|u\|_{H^1} \leq C (\|G_b u\| + \|u\|) \leq C( 1 + |z|) \Bigl( \|G_b - z) u\| + \|u \one_{|x| < 2R}\| + \|u \one_{|x| \geq 2R}\|\Bigr).
$$
Estimate $\|\one_{|x| \geq 2R} u\|$
as follows. 
Choose $\varphi \in C^{\infty} (\R^n)$ 
so that $\varphi = 0$ for $|x| \leq R,$  and $\varphi = 1$ for $|x| \geq 2R.$ 
Define $f:=(G_b - z)u$.
Then $G_b = G_0$ on the support of $\varphi$ and 
$$(G_0 - z) (\varphi u) = \varphi f + [G, \varphi] u.$$
The operator $G_0$ is skew self-adjoint with a dense domain in $L^2(\R^n: \C^r)$ and $\|(G_0 - z)^{-1}\| \leq \frac{1}{|\Re z|}$ for $\Re z \neq 0.$
Thus one obtains
$$
\|\varphi u \| \leq \frac{1}{|\Re z|} \Bigl( \|\varphi f \| + \|[ G, \varphi] u \|\Bigr) \leq \frac{C_1}{|\Re z|} (\|f\| + \|{\bf 1}_{|x| \leq 2R} u \|).
$$
On the other hand, $\one_{|x| \geq 2R} u = \one_{|x| \geq  2R} \varphi u $, 
and this yields the estimate. \hfill\qed

\begin{corollary} For fixed $z \in \C$ with $\Re z < 0$,
$\dim \Ker (G_b - z) < \infty$ and the range of the operator $(G_b - z):\:\hc_0 \cap D(G_b) \longrightarrow \hc_0$ is closed.

\end{corollary}
{\it Proof.} Suppose $\Re z < 0$ and define $B_0 := \{ x \in \Ker (G_b - z):\: \|x\| = 1\}.$ 
Then $B_0 \subset \hc_0.$  Proposition 2.1 implies  that $B_0$ is compact. Indeed, let $\{u_j\}_{j\in \N}$ be  sequence
such that $u_j \in B_0, \: j \in \N.$ Then, from (\ref{eq:4}) one obtains $\|u_j\|_{H^1(\Omega)} \leq C(z),\: \forall j \in \N$. By Rellich theorem we can find a subsequence $\{u_{j_k}\}$ converging in $L^2(\{|x| \leq 2R\})$. The estimate (\ref{eq:4}) shows that
$\{u_{j_k}\}$ is converge in $H^1(\Omega)$. Consequently,  $B_0$ is a compact set and $\Ker (G_b - z)$ is finite dimensional. 

Write $\hc_0 = \Ker (G_b - z) \oplus_\perp B_1$, where $B_1 = (\Ker(G_b - z))^{\perp}$ is the orthogonal complement.  We claim that with a constant $c > 0,$ 
\begin{equation} \label{eq:5}
\|u\| \leq c \|(G_b- z)u\|,\qquad u \in B_1 \cap D(G_b).
\end{equation}
If this inequality were not true, there
would  exist a sequence $v_j \in B_1 \cap D(G_b),\: \|v_j\| = 1$ with
 $(G_b - z)v_j \to 0.$ 
 
 Applying (2.1), and repeating the argument of the proof of Proposition 2.1, it follows  that there exists a subsequence $v_{j_k}$  converging to $w \in B_1$ with
 $\|w\|=1$. Since $(G_b - z)$ is 
 a closed operator, $(G_b - z)w = 0$. 
 Therefore $w \in \Ker(G_b - z) \cap B_1=\{0\}$, so  $w = 0$.
 This is a contradiction and the claim is established. 
 
 To prove that the range is closed, suppose that  $(G_b - z) u_j \to g$ with $u_j = w_j + v_j,\: w_j \in \Ker (G_b - z), \: v_j \in B_1 \cap D(G_b).$ Obviously, $(G_b - z) v_j \to g$ and by (\ref{eq:5}) the sequence $\|v_j\|$ is a Cauchy sequence. Then by (\ref{eq:4}) there exists a convergent subsequence $v_{j_k} \to y$.  Since $(G_b - z)$ is a closed operator, 
 this implies $(G_b - z)y = g$ and the range of $(G_b - z)$ is closed. 
  \hfill\qed

\vskip.2cm

The same argument works for  $\dim \Ker (G_b^* - z)$ and ${\rm Range}\: (G_b^* - z)$.
In particular,   the codimension of ${\rm Range}\: (G_b- z)$ is finite. Consequently, if $z \in \sigma(G_b), \:\Re z < 0,$ and $z$ is not an eigenvalues of $G_d$, 
it follows that ${\rm Range}\: (G_b - z) \not= \hc.$ This means that $z$ is in the residual spectrum of $G_b$ and 
therefore $\bar{z}$ is an eigenvalue of $G_b^*.$

Recall the definition of Fredholm  operators.
 
\begin{definition}  
A bounded operator $L:X \longrightarrow  Y$ from a
Banach space $X$ to Banach space $Y$ is called Fredholm if its kernel is
finite dimensional and its range has finite 
codimension.   The index of a Fredholm operator $L$ is
defined as the difference ${\rm codim}\: {\rm Range}\:\,(L)
-{\rm dim\,Ker}\,(L)$.
\end{definition}

 In the following $\caH_j$ denote  Hilbert spaces and
${\mathcal O}$ a connected open subset of $\C$. 
Recall that a meromorphic operator valued  function $L(\tau) :\caH_1 \to \caH_2, \: \tau \in {\mathcal O}$  is called
finitely-meromorphic in ${\mathcal O}$, if the principal part of the Laurent expansion of $L(\tau)$ at each pole $\tau = \tau_0$
is an operator of a finite rank, that is the coefficients for negative powers of $(\tau- \tau_0)$ are finite-dimensional operators. 

The next analytic Fredholm theorem is a partial case of general results proved in \cite{B}, \cite{K},  where a finitely-meromorphic function $L(\tau): \caH_1 \to \caH_2$ of Fredholm operators $L(\tau)$ is considered.

\begin{theorem} 
Suppose that $L(\tau):\:\caH_1\to \caH_2$ is an analytic function of 
$\tau\in  \caO$ with values in the space of bounded operators from 
$\caH_1$ to $\caH_2$.  Suppose that for all $\tau\in \caO$, the operator $L(\tau)$
is Fredholm.  
Then the set
$$
\big\{
\tau\in \caO\,:\,
L(\tau) {\rm \ \ is \ invertible} \big\}
$$
is either empty or an open subset with discrete
complement.  In the latter case the operator valued function $(L(\tau))^{-1}$
is finitely-meromorphic on $\caO$.
\end{theorem}

\begin{example}   Take
$\caH_1$ to be the set $(H^1(\Omega))^r$ of vector valued 
functions satisfying the dissipative boundary
condition $ f\vert_{\partial \Omega} \in \cn(x)$ 
and belonging to $\caH_0$.
Let  $\caH_2$ be the set of vector valued
$(L^2(\Omega))^r$ functions. When hypotheses $(H)$ and $(H)^*$ 
are satisfied,
$L(\tau):=G_b - \tau:\caH_1\to \caH_2$ 
is an analytic family of bounded operators
and the graph norm of $G_b$ is equivalent to the $H^1(\Omega)$
norm.
 The fact that
 $G_b -\tau:\:\caH_1 \longrightarrow  \caH_2$ is Fredholm
follows from Corollary 2.2. In  Section 4 we will prove that
the  resolvent set of $G_b$ in $\caO:=\{ \tau \in \C: \: \Re \tau < 0\}$ in not empty.
This implies that  
$G_b - \tau$ is invertible on the complement of  
a discrete subset of finite multiplicity eigenvalues in  $\caO$.
\end{example}

\begin{remark}  The Fredholm property of $G_b - \tau$ fails on the imaginary axis,
as it does for the unperturbed operator $G - \tau.$
\end{remark}

\section{Scattering matrix $\Sc(z)$, invertibility of $\Sc(0)$}

Introduce  the wave operators related to $U_0(t)$ and $V(t),$
assuming  that $(H)$ and $(H)^*$ are satisfied. Consider the operator $J:\: \hc \longrightarrow H_0$ extending $f \in \hc$ as 0 for $x \in K$ and let $J^*:\: H_0 \longrightarrow \hc$ be the adjoint of $J.$ Let $\Ha = (\Ker G_0)^{\perp}.$ The wave operators are defined by
 $$
W_{-}f = \lim_{t \to +\infty}  V(t) J^* U_0(-t) f
\qquad
W_{+} f = \lim_{t \to +\infty} V^*(t) J^* U_0(t) f,\qquad f \in \Ha.
$$
It is not difficult to prove the existence of $W_{\pm}$ (see, for instance, \cite{GS} and Chapter III in \cite{P2}). 
In addition,
$$
 \overline{{\rm Ran}\: W_{\pm}}\  \subset \ 
 \Hp\,.
 $$
It was proved in \cite{P1} that if $G_b$ has eigenvalue in $\Re z < 0$ with eigenfunction $f$ whose  support is not compact, then the wave operators are not complete and ${\rm Ran}\: W_{-} \not= {\rm Ran}\: W_{+}.$ 
In the next section we prove that the eigenfunctions with eigenvalues in $\Re z < 0$ never have compact support.

 The corresponding operators with the roles
 of free and perturbed interchanged  
 are  
$$
Wf := \lim_{t \to \infty}  U_0(-t) J V(t) f,\qquad
W_1f := \lim_{t \to \infty}  U_0(t) J V^*(t) f,
\qquad f\in \Hp\,.
 $$
 In \cite{GS}, \cite{P2},
 the existence of these limits is proved 
assuming the hypothesis $(H)$ , $(H)^*$.
 
The Hilbert spaces $\Ha$ and $\Hp$ are invariant 
under the action of $U_0(t)$ and $V(t)$ respectively.
This is the setting of 
two spaces scattering theory with
scattering operator defined 
as   $S := W \circ W_-$. We have $\overline{{\rm Ran}\: W} \subset \Ha$ and $S$ is a bounded operator from $\Ha$ to $\Ha$.

We recall some elements of   Lax-Phillips scattering theory
that are needed.
The translation representation 
$${\mathcal R}_n: \Ha \longrightarrow (L^2(\R \times \sn))^d$$
 of $U_0(t),$ 
 involves both the Radon transform 
 $(Rf)(s, \omega)$ of $f(x)$  and eignevalue-eigenvector pairs
 $\tau_j(\xi),r_j(\xi)$.
 With Rank $A(\xi) = r - d_0 = 2d > 0$ for $\xi \not= 0$,
 define $\tau_j(\xi),\: j = 1,...,2d,$ to 
 be the non-vanishing eigenvalues of $A(-\xi)$ in decreasing order,
$$
\tau_1(\xi) \geq ...\geq \tau_d(\xi) > 0 > \tau_{d+1}(\xi) \geq ...\geq \tau_{2d}(\xi),\quad \xi \not= 0.
$$
Choose measurable normalized eigenvectors  $r_j(\xi), \: j = 1,...,2d,$  of $A(-\xi)$ with eigenvalues  $\tau_j(\xi)$. Then $\Rc_n$ has the form (see Chapter VI in \cite{LP1} for $d_0 = 0$ and Chapter IV in \cite{P2} for $d_0 > 0$)
$$(\Rc_n f)(s, \omega) = \sum_{j=1}^d \tilde{k}_j(s, \omega) r_j(\omega ),\qquad
\tilde{k}_j(s, \omega) := \tau_j(\omega)^{1/2} k_{j}(s \tau_j(\omega), \omega),\: j = 1,...,d,$$
where 
$$
k_j(s, \omega)\ :=\
 2^{-(n-1)/2}D_s^{(n-1)/2} \langle (Rf)(s, \omega), r_j(\omega) \rangle
\,.
$$  
 $\Rc_n$ is an isometry $\Ha\to (L^2(\R \times \sn))^d$, and $\Rc_n U_0(t) = T_t \Rc_n,\: \forall t \in \R$, where $T_t g = g(s- t, \omega).$

\begin{definition} 
Define the slowest nonzero sound speed $v_{min} := \min_{|\omega| = 1} \tau_d (\omega) > 0.$ 
The set 
$D_{\pm}\subset H_0$ is the set of  $f$ satisfying 
$$
U_0(t)f = 0 \quad {\rm for}\quad |x|< \pm v_{min} t,\quad \:\pm t > 0.
$$
Equivalently $f \in D_{\pm}$ if and only if $\Rc_n f(s, \omega) = 0$ for $\mp s > 0.$
Define
$$
D_{\pm}^{a}
\  :=\
 U_0(\pm a/v_{min}) D_{\pm},\quad a \geq \rho.
$$
\end{definition}

Then  $D_{\pm}^{\rho} \subset \Hp$ (see Lemma 4.1.6 in \cite{P2}),
$$
J(\overline{D_{+}^{\rho} \oplus D_{-}^{\rho}}) 
\ =\
 \overline{D_{+}^{\rho} \oplus D_{-}^{\rho}},$$
$$ J^*(\overline{D_{+}^{\rho} \oplus D_{-}^{\rho}}) 
\ =\
 \overline{D_{+}^{\rho} \oplus D_{-}^{\rho}}\,.
$$
It follows that  the  operators $W_{\pm}$ and  $W$ defined above coincide with the Lax-Phillips wave operators  \cite{LP2} related to $\Ha$ and $\Hp$.

Consider $D_{\pm}^{a},\: a \geq \rho,$ as subspaces of $\Hp$ and 
 introduce the orthogonal projectors $P_{\pm}^{a}$ on the orthogonal complements of $D_{\pm}^{a}$ in $\Hp.$ It is easy to see (cf. Chapter IV, \cite{P2}) that 
$$
(i) \  V(t) D_{+}^{\rho} \subset D_{+}^{\rho} \quad
{\rm and}\quad
V^*(t) D_{-}^{\rho} \subset D_{-}^{\rho},\: t \geq 0,
$$
$$
(ii)\ \bigcap_{t \geq 0} V(t) D_{+}^{\rho} = \{0\}
\quad
{\rm and}\quad
\bigcap_{t \geq 0}  V^*(t) D_{-}^{\rho} = \{0\}\,.$$
Next we prove
$$
(iii)\:  \lim_{t \to \infty} P_{+}^{\rho} V(t) f = 0
\quad
\forall f \in \Hp
\quad
{\rm and}\quad\lim_{t \to \infty} P_{-}^{\rho} V^*(t)f = 0,\quad \forall f \in \Hp.
$$

The existence of the operator $W$ implies the first assertion.  
 In fact, given $g \in \Hp$, there exists $f \in H_0^{ac}$ such that
$$\lim_{t \to \infty} \|V(t) g - (U_0(t) f)\vert_{\Hp}\|_{H} = 0.$$
Thus it suffices to prove that
 $\lim_{t \to \infty} \|P_{+}^{\rho} (U_0(t) f) \vert_{\Hp})\|_{H} = 0.$ It is
 clear that this  holds for $f = U_0(s) h,\: h \in D_{+}^{\rho}.$   Since 
$$
\overline{ \bigcup_{t \in \R} U_0(t) D_{+}^{\rho}} = H_0^{ac}\,,
$$
the assertion follows 
for every $f \in H_0^{ac}$.
Thus $\lim_{t \to \infty} P_{+}^{\rho} V(t) g = 0.$ 

For the other relation in (iii),
 use the existence of the operator $W_1$.

 Properties  (i)-(iii), justify the application of  
  the abstract Lax-Phillips scattering theory for dissipative operators \cite{LP2}.

Using the translation representation $\Rc_n$ of $U_0(t)$, consider the scattering operator
$$
\tilde{S}: (L^2(\R \times \sn))^d \rightarrow (L^2(\R \times \sn))^d,
\qquad
\tilde{S} := \Rc_n S \Rc_n^{-1}\,.
$$
The kernel of $\tilde{S} - I$ is a matrix-valued distribution $K^{\#}(s- s', \theta, \omega)$ and 
$$
K^{\#}(s, \theta, \omega) 
\ =\
 \Bigl(S^{ji}(s, \theta, \omega)\Bigr)_{j,i = 1}^d
$$
 is called {\it scattering kernel}. The relation $(S - I)(D_{-}^{\rho}) \subset (D_{+}^{\rho})^{\perp}$ (see Lemma 4.1.7 in \cite{P2}) implies 
$$
K^{\#}(s, \theta, \omega)  = 0 \quad {\rm for}\quad s > 2\rho/v_{min},
$$
and the Fourier transform  $\hat{K}^{\#}(z, \theta, \omega)$ of $K^{\#}(s, \theta, \omega)$ with respect to $s$ is analytic for $\Im z < 0.$ The same is true for the operator-valued function $\Sc(z) - I$ with kernel $\hat{K}^{\#}(z, \theta, \omega).$ The operator $\Sc(z)$ is called the {\it scattering matrix} (cf. for instance \cite{LP2}). 
More precise information about the support of $K^{\#}(s, \theta, \omega)$ and the maximal singularity of $K^{\#}(s, -\omega, \omega)$ is given in \cite{MT}, \cite{P1} for the case $d_0 = 0$ and in  \cite{P3} for the case studied in this paper. In particular, in  \cite {P3}  one proves a representation of the scattering kernel and taking the Fourier transform with respect to $s$ 
it follows  that  $\Sc(z) = I + \Kc(z)$ for $\Im z < 0$ is an analytic operator-valued function and $\Kc(z)$ is an operator-valued function with values in the space of Hilbert Schmidt operators in $(L^2( \sn))^d.$

  To obtain a representation of $\Kc(z)$,
  recall Theorem 15.5 in \cite{P3}. This result was proved 
  when $A(\xi)$ has characteristic roots with constant multiplicities for $\xi \not= 0$.
  The proof works without any changes when
   only zero is a characteristic root of constant multiplicity. Let  
$$
w_k^o (t, x, \omega) = \tau_k(\omega)^{1/2} \delta^{(n-1)/2}( \langle x, \omega \rangle - \tau_k(\omega)t) r_k(\omega),\: k = 1,...,d,$$
and consider the solution $w_k^s(t, x, \omega)$ of the problem
\begin{equation} \label{eq:3.1}
\begin{cases} (\partial_t - G)w_k^s  = 0 \: {\rm in} \: \R \times \Omega,\\
w_k^s + w_k^o   \in {\mathcal N}(x) \: {\rm on}\: \R \times \partial \Omega,\\
w_k^s \big\vert_{t \leq - \frac{\rho}{v_{min}}} = 0.\end{cases}
\end{equation}
 The scattering kernel $K^{\#}(s, \theta, \omega)$ computed with respect to the basis $\{r_j(\omega)\}_{j=1}^d$ has matrix elements
 $S^{j k}(s, \theta, \omega)$ equal to
\begin{equation} 
\label{eq:3.2}
d_n^2 \tau_j(\theta)^{1/2}
\hskip-4pt
\int\limits_{\R \times \partial \Omega}
\hskip-2pt \delta^{(n-1)/2}(\langle x, \theta \rangle - \tau_j(\theta)(s + t))\Big\langle r_j(\theta), A(\nu(x))  (w_k^o + w_k^s)(t, x, \omega)\Big\rangle dt dS_x,
\end{equation}
where $d_n \neq 0$ is a constant depending only on $n$ and the integral is taken in sense of the distributions.  Thus for $\Im z < 0$, taking the Fourier transform with respect to $s$, one obtains
\begin{eqnarray} \label{eq:3.3}
\hat{S}^{ji} (z, \theta, \omega) = (-\ii z)^{(n-1)/2} d_n^2 \tau_j(\theta)^{1/2} \\
\times \int\limits_{\R \times \partial \Omega} \exp\Bigl( \ii z \Bigl(\frac{\langle x, \theta \rangle}{\tau_j(\theta)} - t \Bigr)\Bigr)  \Big\langle r_j(\theta), A(\nu(x))  (w_k^o + w_k^s)(t, x, \omega)\Big\rangle dt dS_x. \nonumber
\end{eqnarray}

The above equality for $\Im z < 0$  yields
\begin{equation} 
\label{eq:3.4}
\hat{K}^{\#}(z, \theta, \omega) = c_n z^{(n-1)/2} K_1(z, \theta, \omega), \qquad c_n \not= 0
\end{equation}
with 
$K_1(z, \theta, \omega)$ an analytic matrix-valued function for $\Im z < 0.$\\

 Since $w_k^o(t, x, \omega)\big\vert_{\R \times \partial \Omega}$ has compact support with respect to $t$, the integral
$$\int\limits_{\R \times \partial \Omega} \exp\Bigl( \ii z \Bigl(\frac{\langle x, \theta \rangle}{\tau_j(\theta)} - t \Bigr)\Bigr)
\Big\langle r_j(\theta), A(\nu(x)) w_k^o(t, x, \omega)\Big\rangle dt dS_x,$$
 yields an analytic function in $\C$. In the integral in (\ref{eq:3.3}) over $\R$  involving $w_k^s$, for $\Im z < 0$ take the inverse Fourier transform ${\mathcal F}_{t \to z}$ with respect to $t$ in the sense of distributions and denote
$$v_k^s (z, x, \omega) = {\mathcal F}_{t \to z} w_k^s(t, x, \omega).$$
 Then $v_k^s(z, x, \omega)$ is a solution of the problem
$$ \begin{cases} (\ii z - G) v_k^s = 0 \: {\rm in} \: \Omega,\\
   \Bigl(v_k^s  + (\ii z)^{(n-1)/2} e^{-\ii z \frac{\langle x, \omega \rangle}{\tau_k(\omega)}} r_k(\omega)\Bigr)\bigg\vert_{\partial \Omega} \in {\mathcal N}(x)\: {\rm on}\: \partial\Omega. \end{cases}$$

Let $\varphi(x) \in C_0^{\infty}(\R^n)$ be  a function such that $\varphi(x) = 1$ for $|x| \leq \rho,\: \varphi(x) = 0$ for $|x| \geq 2 \rho.$ Write $v_k^s(z, x, \omega) = V_k^s(z, x, \omega) - (\ii z)^{(n-1)/2}\varphi(x)  e^{-\ii z \frac{\langle x, \omega \rangle}{\tau_k(\omega)}} r_k(\omega)$ and deduce
$$\begin{cases} (\ii z - G) V_k^s = -[G, \varphi] (\ii z) ^{(n-1)/2} e^{-\ii z \frac{\langle x, \omega \rangle}{\tau_k(\omega)}} r_k(\omega)\: {\rm in} \: \Omega,\\
V_k^s \in {\mathcal N}(x)\: {\rm on}\: \partial \Omega. \end{cases}$$

Consequently, setting $-[G, \varphi] = \psi(x)$, one has
$$V_k^s (z, x, \omega) = (\ii z- G_b)^{-1} \psi(x)  (\ii z) ^{(n-1)/2} e^{-\ii z \frac{\langle x, \omega \rangle}{\tau_k(\omega)}} r_k(\omega)$$
and
$$V_k^s (z, x, \omega)\big\vert_{\partial \Omega}  = (\ii z)^{(n-1)/2} \Bigl[\varphi(x)(\ii z- G_b)^{-1} \psi(x) e^{-\ii z \frac{\langle x, \omega \rangle}{\tau_k(\omega)}} r_k(\omega)\Bigr] \bigg\vert_{\partial \Omega}.$$
In conclusion, the kernel $K_1(z, \theta, \omega)$ for $\Im z < 0$ is given by a sum of a function analytic with respect to $z$ in $\C$ and an integral over $\partial \Omega$ involving the cut-off resolvent $\varphi(x) (\ii z - G_b)^{-1} \psi(x).$ A similar argument has been used for the scattering amplitude $s(z, \theta, \omega)$ related to the wave equation in \cite{PS}.\\

To introduce the scattering resonances, 
 consider in $\Hp$ the semigroup  of contractions $Z(t)$ in
 the space 
 $K^{a} := \Hp \ominus (D_{-}^a \oplus D_{-}^a), \: a \geq \rho$,
 $$
 Z_a(t)\ :=\
  P_{+}^{a} V(t) P_{-}^{a} \ :=\
   e^{tB_a},
  \qquad
  t \geq 0.
  $$ 
  Exploiting the condition $(H)$, as in \cite{LP1}, \cite{LP2},
 it follows  that the spectrum 
 of $B_a$ is discrete and formed  only by isolated eigenvalues of finite multiplicity in $\Re z < 0$ and this spectrum is independent on the choice of $a$.  Corollary 4.8 in \cite{LP2} 
  implies that {\it  the resolvent of $B_a$ is a meromorphic function in $\C$.} Now choose $a \geq 2 \rho.$ Then $\varphi P_{+}^a = \varphi,\: 
P_{-}^a \psi(x) = \psi(x)$ and for $\Im z < 0,$
$$\varphi (\ii z - B_a)^{-1} \psi = \int_0^{\infty} \varphi P_{+}^a V(t) P_{-}^a \psi dt = \varphi(\ii z - G_b)^{-1} \psi.$$
The left hand side is analytic for $\Im z \leq 0$, so $\varphi (\ii z - G_b)^{-1} \psi$ is analytic for $\Im z \leq 0$. 
This implies that $\lim_{z \to 0, \: \Im z < 0} K_1(z, \theta, \omega) = K_1(0, \theta, \omega)$ is  a kernel of a bounded operator in $(L^2(\sn))^d$ and ${\mathcal S}(0)$ is invertible.

\section{The spectrum of $G_b$}

 {\bf Fix} $a =2\rho$ and write simply $Z(t),\: B$ instead of $Z_a(t),\: B_a$. 
 
Passing to the translation representation,
yields the decomposition $\Hp = K^{a} \oplus D_{-}^a \oplus D_{+}^a.$

\begin{definition}
A function $f$ is {\bf outgoing (resp. incoming)}
 if the component of $f$ in $D_{-}^a$ (resp. $D_{+}^a$) vanishes. 
 \end{definition}
 
 In \cite{LP2} the incoming (outgoing) functions are defined with respect to the components in $D_{+}^{\rho}$ ($D_{-}^{\rho}$). For our analysis the role of $\rho$ is played by $a = 2\rho$. To apply the results of Section 5 in \cite{LP2}, we need to prove that neither $G_b$ nor $G_b^*$ have eigenfunctions in $\Hp$ 
that are both incoming and outgoing.
 A stronger result is Corollary \ref{cor:noeigen}
below. First we prove the following

\begin{theorem} Let  $f \in H_b^{\perp}.$ Then the following conditions are equivalent.

\vskip.1cm
$(a)$ there exists $b \geq 2\rho$ such that $f \perp D_{-}^b$ and $\lim_{t \to \infty} V(t) f  = 0.$
\vskip.1cm
$(b)$ $V(t)f$ is a disappearing solution.
\end{theorem} 

\noindent
This results is similar to Theorem 1 of Georgiev \cite{G2} (see also Theorem 4.3.2 in \cite{P2}) established in the case $\det A(\xi) \not= 0$ for $\xi \not = 0.$ For non elliptic symmetric systems
new ideas are needed.

{\it Proof of Theorem $4.2$.} The implication (b) $\Rightarrow$ (a) follows  from the finite speed of propagation. We prove that (a) $\Rightarrow$ (b).  Assume (a). First treat the case when $f \in \cap_{j=1}^{\infty} (G_b)^j \cap \Hp \cap (D_{-}^b)^{\perp}.$ Set $u(t, x) = V(t)f$ and consider the solution
$$
\tilde{u}(t, x) 
\ :=\
 V(t) G_bf\ =\
  G_b V(t) f
  \  =\
   \partial_t (V(t)f).
$$
  Since $f \perp D_{-}^{\rho},$ we have $V(t) f \perp D_{-}^{\rho},\: t \geq 0.$ Taking the derivative with respect to $t$, yields
   $\tilde{u}(t, x) \perp D_{-}^{\rho}.$ 
   
   We claim that $\lim_{t \to + \infty} \tilde{u}(t, x) = 0$.
    Our hypothesis implies
 $$
 \|\partial_t^j (V(t)f) \|
 \ =\
  \|V(t) (G_b)^j f\|
  \  \leq\
   C_j,\qquad\forall t \geq 0, \ \forall j \in \N.$$
 To prove the claim, suppose that there exists a sequence
$t_k \to + \infty$ such that $\|u_t(t_k, x)\| \geq \delta > 0, \: \forall k \in \N.$ Then for $\xi \leq \eta \leq t_k$ we get
$$
\|u_t(t_k, x) - u_t(\eta, x) \| 
\ =\
 \big \| \int_{\eta}^{t_k} u_{tt}(y, x) dy \big \|
 \  \leq\
  C_2 (t_k - \xi) 
  $$
   and, provided $0 \leq t_k - \xi \leq \frac{\delta}{2 C_2},$ one has
$$
\|u(t_k, x) -u(\xi, x)\| 
\ =\
 \big\|\int_{\xi}^{t_k} u_t(y, x) dy \big \| 
 \ \geq\
  (t_k - \xi) \delta - (t_k- \xi)^2 C_2 
  \ \geq\
   (t_k- \xi) \frac{\delta}{2}.
  $$
 This contradicts the hypothesis $u(t_k, x) \to 0$ and the claim is proved.

Choose a function
 $\varphi(x) \in C^{\infty}(\R^n)$ with
  $\varphi(x) = 1$ for $|x| \geq 2\rho,\: \varphi(x) = 0$ for $|x| \leq \rho.$ 
  With  $v(t,x) := \varphi(x) \tilde{u}(t, x) $ compute
$$ 
v(t, x) 
\ =\
 \varphi V(t) G_b f 
 \ =\
  G \varphi V (t) f + [\varphi, G] V(t) f 
\ :=\
 w(t, x) + [\varphi, G] V(t) f,
$$
the last equality defining $w$.

Next we will
prove that $w(t, x) = 0$ for $|x| \geq b,\: t \geq 0.$ 
Start with  $G_b \varphi = G \varphi$. 
It is clear that
$w(t, x) \in  (\Ker G)^{\perp} = \Ha$, so we may consider the translation representation $\Rc_n( w(t, x)) = m(t, s, \omega)$ of $w(t, x).$  
Since $f \in H_b^{\perp},$ there exists a sequence $t_k \to \infty$ such that $\lim_{t_k \to \infty} \|V(t_k) f\|_{L^2(|x| \leq 2\rho)} = 0$  (see Proposition 3.1.9 in \cite{P2}). On the other hand, $\|\tilde{u}(t_k, x)\|_{L^2(\Omega)} \to 0$ as $t_k \to \infty$, so
 $\lim_{t_k \to \infty} \|w(t_k, x)\|_{L^2(\Omega)} = 0$.

Write, with the last equality defining $g$,
$$
(\partial_t - G) w(t,x) 
\ =\
  -\sum_{j = 1}^n (A_j \varphi_{x_j}) V(t)G_b f + [G, [\varphi, G]] V(t)f 
\ := \
g(t,x)\,.
$$
 Then $g(t, x) \in (L^2(\R^{+} \times \R^n))^r$ vanishes for $|x| \geq 2 \rho.$
Applying the transformation $\Rc_n$ to both sides of the above equality and setting 
$$
m_j(t, s, \omega) = \langle m(t, s, \omega), r_j(\omega) \rangle,
\quad
 l_j(t, s, \omega) = \langle \Rc_n(g)(t, s, \omega), r_j(\omega) \rangle, \ \
  j = 1,...,d,$$
yields 
$$
(\partial_t + \tau_j(\omega) \partial_s) m_j(t, s, \omega) 
\ =\
 l_j(t, s, \omega),\quad j =1,...,d.
 $$
 Fix $t_1 \geq 0$ and $1\le j\le  d$.
  Since $l_j(t, s, \omega) = 0$ for $|x| \geq 2 \rho,\:t \geq t_1$ and $s \geq 2 \rho + \tau_j(\omega) (t- t_1)$,
  the above equations yield
$$
m_j(t, s, \omega) \ =\
 m_j(t_1, s - \tau_j(\omega) (t-t_1) , \omega)\,.
 $$
 This implies
$$\int_{\R} \int_{\sn}  |m_j(t, s, \omega)|^2 ds d\omega \geq  \int_{2\rho}^{\infty} \int_{\sn}|m_j(t_1, s, \omega)|^2 ds d\omega.$$

Letting $t_k \to \infty$, one gets
 $m_j(t, s, \omega) = 0$ for $s \geq 2\rho.$ 
 Since 
  $V(t)f \perp D_{-}^b$ 
  and 
  ${\rm supp}\,(1- \varphi(x))\subset\{|x| \leq 2\rho\},$
  this implies that 
   $w(t, x) \perp D_{-}^b$ for $t \geq 0$.
   Therefore, $m(t, s, \omega) = \Rc_n (w(t, x)) = 0$ for $s \leq -b.$ 
   
   Next repeat the argument of Lemma 2.2 in \cite{G2} based on the following
\begin{lemma}
 [\cite{G1}] Let $F \in \bigcap_{j=1}^{\infty} D(G^j) \cap \Ha$ and let $\Rc_n F = 0$ 
 for $|s| \geq b.$ Then $F = 0$ for $|x| \geq b$ if and only if
$$
\int\int_{\R \times \sn} [A(\omega)]^k 
\big [ (\Rc_n F)(s, \omega) + (-1)^{(n-1)/2} (\Rc_n F)(-s, -\omega)\big ] 
s^a Y_m(\omega) ds d\omega = 0
$$
for $a = 0,1,2,...$ and any spherical harmonic function $Y_m(\omega)$ 
of order $m \geq a + k +(3-n)/2.$
\end{lemma}

We conclude that $w(t, x)$ and  $v(t, x)$ have support in $\{x: |x| \leq b\}$ for all $t \geq 0.$ Thus $\tilde{u}(t, x) = \partial_t (V(t) f) = 0$ for $|x| \geq b$. Consequently, our hypothesis implies that $u(t, x) = 0$ for $|x| \geq b.$

The next step is to show that  $\supp\:\: u(t, x) \subset \{x: |x| \geq b\}$ 
implies that $u$ is  a disappearing solution. We know that $\partial_t \tilde{u} - G \tilde{u} = 0$. As we have mentioned in Introduction, according to \cite{M}, under the hypothesis that $A(\xi)$ has constant rank for $\xi \neq 0$, it is possible to show that  there is an $(r \times r)$ matrix-valued 
polynomial $Q(\xi)$ such that 
\begin{equation} \label{eq:4.1}
\Ker Q(\xi) = {\rm Image}\: A(\xi),\quad \xi \neq 0.
\end{equation}
Since there is no proof of this property in \cite{M}  and since we wish to exploit the structure of $Q(\xi)$, for the sake of completeness
 we present a proof of (\ref{eq:4.1}). Recall that $r - d_0 = 2d > 0$. For $\xi \neq 0$ let
$$
P(z, \xi) = R(z, \xi)z^{d_0} 
\ =\
  \prod_{j=1}^{2d} (z + \tau_j (\xi)) z^{d_0}
  $$
be the characteristic polynomial of $A(\xi),$  where $-\tau_j(\xi)$ are the non-vanishing eigenvalues of $A(\xi)$ repeated with their multiplicities. Define $Q(\xi) = R(A(\xi), \xi),\: \xi \neq 0.$ 

Next show that $Q(\xi)$ is a polynomial $\sum_{|\alpha| \leq 2d} B_{\alpha}\xi^{\alpha}$ with matrix coefficients $B_{\alpha}.$  Since $A^j(\xi)$ are polynomials with matrix coefficients, it is sufficient to prove that 
$R(z, \xi) = \sum_{j=0}^{2d} c_j(\xi) z^j$ has coefficients $c_j(\xi)$ 
that are polynomials. This follows from
$$
\det (zI - A(\xi)) 
\ =\
 \sum_{j=0}^{2d} c_j(\xi) z^{j + d_0}
 $$
 upon comparing the coefficients of $z^{j+d_0}$ in both sides.
The  Cayley-Hamilton theorem implies that
$$
P(A(\xi), \xi) \ =\
  Q(\xi) A^{d_0} (\xi)
  \  =\ 
   0.
$$
Passing to a diagonal form of $A(\xi)$, 
shows that ${\rm Ran}\: A(\xi) = {\rm Ran}\:A^{d_0} (\xi).$ Thus ${\rm Ran}\: A(\xi) \subset \Ker Q(\xi),\:\xi \neq 0.$ To establish the
 opposite
  inclusion, assume that $h \in \Ker Q(\xi)$ and write $h = h_1 + h_2$ with $h_1 \in \Ker A(\xi),\: h_2 \in {\rm Ran} \: A(\xi).$ Then $Q(\xi) h = \prod_{j=1}^{2d}  \tau_j(\xi) h_1 = 0$ and we conclude that $h_1 = 0.$\\

  Let $Q(D_x)$ be the operator with symbol $Q(\xi)$ and let $L = (\partial_t -G)^{4d} + Q^2.$ Then $Q \tilde{u} = Q G V(t)f = 0$ and we get $L \tilde{u}= 0 $.  The symbol of the operator $L$ with constant coefficients is 
$$ (\tau I- A(\xi))^{4d} +  Q^2(\xi) =  l(\tau, \xi).$$
First, show that $\det l(0, \xi) \neq 0$ for $\xi \neq 0.$ In fact, if for $\xi_0 \neq 0$ there exists a vector $v \neq 0$ such that $A^{4d}(\xi_0)v + Q^2(\xi_0) v = 0$, taking the scalar product by $v$,
yields $A(\xi_0) v = Q(\xi_0) v = 0.$ This implies $v \in \Ker A(\xi_0) \cap {\rm Ran} \: A(\xi_0)$ and 
therefore $v = 0$ which is a contradiction. Next, if for $\tau \neq 0$ and $\xi \neq 0$ we have $\det l(\tau, \xi) = 0$, by the same argument we deduce that $\det(\tau I - A(\xi)) = 0$ and 
therefore $\tau = -\tau_j(\xi)$ for some
$j = 1,...,2d.$ If $\tau < 0,$ we have $\tau = - \tau_j(\xi)$ for some $j = 1,...,d.$ On the other hand, if $\tau > 0$, we get
$\tau = - \tau_j(\xi) = \tau_{2d - j +1}(-\xi)$ for some $j = d+1,...,2d.$
Thus 
$$
\frac{|\tau|}{|\xi|} 
\  \geq\
 \min \{\tau_j(\omega),\: \omega 
 \ \in\
  \sn, j = 1,...,d\} 
  \ =\  v_{\min} \ >\  0.
  $$

Denote by $S(x, r)$ the ball $\{y \in \R^n: |y - x| \leq r\}$.

\begin{proposition} 
Suppose that 
 $S(z, r_0) \subset \Omega$,  $|z - z_1| = \12 r_0 > 0$
 with $0 < r_1 < r_0.$ Assume that $U \in {\mathcal D}'(\R_t \times\Omega)$ is a solution of $L U = 0$ such that $U(t, x) = 0$ for $t \geq t_0 \geq 0,\: x \in S(z, \frac{1}{2} r_0).$ Then
\begin{equation}
U(t, x) = 0 \: {\rm for}\ \: x \in S(z_1, \frac{1}{2} r_1),
\qquad
 t \geq t_0 + \frac{1}{2 v_{\min}} r_0.
\end{equation}

\end{proposition}

{\it Proof.} Let $\Pi$ be a characteristic hyperplane for $L$ with normal $N = (\tau, \xi) \in \R^{n+1} \setminus \{0\}.$ Then 
$\det l(\tau, \xi) = 0$ so
 $|\tau| \geq v_{\min}|\xi|$. Thus, if $|\tau| < v_{\min} |\xi|$ the hyperplane $\Pi$ is not characteristic. Since $r_1 < r_0,$ a simple geometric argument shows that for $s \geq t_0$ every characteristic hyperplane for $L$  which intersects the convex set $\{(t, x):\: |x- z| \leq \12 r_0 + v_{\min} \frac{r_1}{r_0}(t - s), \:s\leq t \leq s + \frac{1}{2 v_{\min}} r_0\}$ intersects also the convex set
$\{ (t, x):\:|x - z| \leq \12 r_0, \:s \leq t \leq s + \frac{1}{2 v_{\min}} r_0\}.$ Thus we can apply F. John's global Holmgren theorem (see for instance, Chapter 1, Corollary 9 in  \cite{Rauch} or Theorem 8.6.8 in \cite{H}) to conclude that $U(t, x) = 0$ for $x \in S(z_1, \12 r_1),\: t \geq t_0 + \frac{1}{2 v_{\min}} r_0.$
For convenience of the reader the non-characteristic deformations of the boundary
  $|x|=\12 r_0, t\ge 0$ are sketched in Figure 1.\hfill\qed

\begin{figure}[ht!]
\begin{center}
\includegraphics[width=5.5cm]{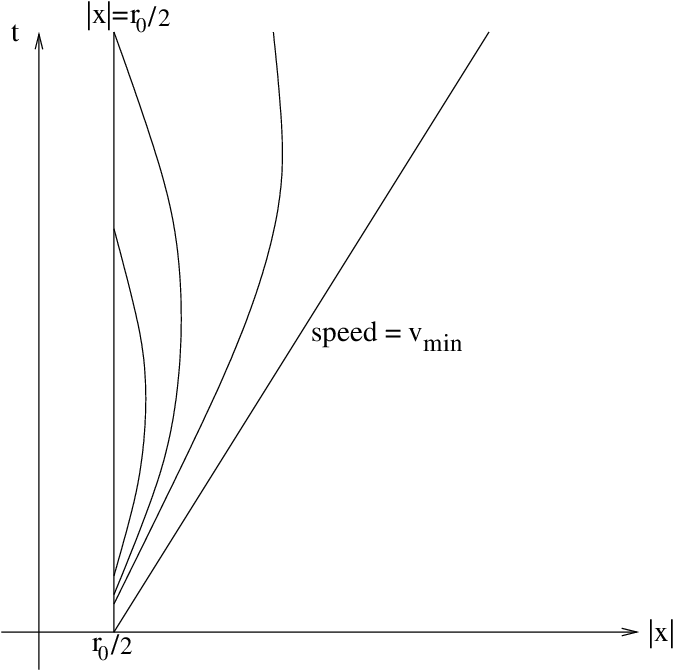}
\caption{Deformation of the boundary}
\end{center}
\label{fig:domain}
\end{figure}

Now suppose that $\tilde{u}(t, x) = 0$ for $|x| \geq b, t \geq 0.$
Fix a point $x_0$ and $r_0$ so that $S(x_0, \12 r_0) \subset \{|x| > b\}$ and $S(x_0, r_0) \subset \Omega.$ Since $\Omega$ is connected, given a point $\hat{x} \in \Omega$, there exists a path $\Gamma_{\hat{x}} \subset \Omega$ with length less than a fixed number $L_0 > 0$ independent on $\hat{x}$ and points $x_1, x_2,..., x_m = \hat{x}$ on the path $\Gamma_{\hat{x}}$ so that
$$S(x_j, r_j) \subset \Omega, \: j = 0, 1,...,m-1,$$
$$r_0 > r_1 >...> r_{m-1} > 0$$
with $r_j = 2|x_{j+1} - x_j|,\: j = 0, 1,..., m-1.$ Assume that
$$
L \tilde{u} = 0\ \  {\rm for}\ \ 
 |x- x_j| \leq  r_j,
\
 t \geq 0,
$$
$$
\tilde{u}(t, x) = 0 \ \ {\rm for} \ \  x \in S(x_j, \12 r_j),\quad
 t \ \geq\  \frac{1}{2 v_{\min}}\,\sum_{k= 0}^{j- 1} r_j.$$

Applying Proposition 4.4 with $U = \tilde{u}$, 
shows that
$$
\tilde{u}(t, x) = 0 
\qquad {\rm for}\qquad   x \in S(x_{j+1}, \12 r_{j+1}),\quad
t \geq 
\frac{1}{2 v_{\min}}\sum_{k= 0}^j r_j.
$$ 
Consequently,
$\tilde{u}(t, \hat{x}) = 0$ for $t \geq \frac{1}{2 v_{\min}} \sum_{k=0}^{m-1} r_k = \frac{1}{2 v_{\min}} L_0$. This argument works for every $x \in \Omega$, hence $\tilde{u}(t, x)$ as well $u(t, x)$ are disappearing. Moreover, the constant $T_0 = \frac{1}{2 v_{\min}} L_0$ depends on $b$ and $\Omega$, so it is independent on $f$.

 Next  treat the general case.
  For every fixed $\epsilon > 0$, we construct a sequence $\{\varphi_{\epsilon}\}$ such that
$$ 
\varphi_{\epsilon} \in \Hp \cap \Bigl( \cap_{j=1}^{\infty} D(G_b^j)\Bigr) \cap D_{-}^{\rho},\qquad
 \lim_{t \to \infty} V(t) \varphi_{\epsilon} = 0
$$
and $\|\varphi_{\epsilon} - f \| < \epsilon.$ The construction of this sequence is given in the proof of Theorem 1 in \cite{G2} and for the reader's
convenience, we sketch it. Assume that a sequence $f_0, f_1,...,f_p,$ is defined so that $f_0 = f$ and  
$$
f_j \in D(G_b^j) \cap \Hp \cap (D_{-}^b)^{\perp},\qquad j = 1,...,p.
$$
Set
$$
f_{p+1} \ :=\
 \frac{1}{\epsilon_{p+1}} \int_0^{\epsilon_{p+1}} V(\tau) f_p
 \ d\tau.
 $$
Then 
$f_{p+1} \in D(G_b^{p+1})$ and 
$$G_b f_{p+1} = \frac{1}{\epsilon_{p+1}}\Bigl( V(\epsilon_{p+1}) f_p - f_p\Bigr).$$
Next 
choose $\epsilon_{p+1} > 0$ so that $\epsilon_{p+1} < \epsilon/2$ and 
$$
\|G_b^k(f_{p+1} - f_p)\| 
\ \leq\
 \frac{\epsilon}{2^{p+1}}\quad
 {\rm for} 
 \quad
  k =0,1,...,p.$$
This is possible since $V(t)$ is strongly continuous at 0 and
$$
\|G_b^k(f_{p+1} - f_p) \| 
\ \leq \
\max_{0 \leq \tau \leq \epsilon_{p+1}} \|V(\tau) G_b^k f_p - G_b^k f_p\|.
$$
Now for fixed integers $N \geq 0,\: p \geq N,\: \mu \geq 0$ we get $\|G_b^N (f_{p+ \mu} - f_{\mu})\| < \epsilon/2^p$ and since the operators $G_b^p$ are closed, we can find $\varphi_{\epsilon} \in \cap_{j=1}^{\infty}  D(G_b^j) \cap \Hp \cap (D_{-}^b)^{\perp}$ so that 
$\lim_{p \to \infty,\: p \geq N} G_b^N f_p = G_b^N \varphi_{\epsilon},\: \forall N \geq 0.$ Finally, $\lim_{p \to \infty} f_p = \varphi_{\epsilon}$, 
implies that
 $\lim_{t \to \infty} V(t) \varphi_{\epsilon} = 0.$

Applying the above argument, we conclude that $V(t) \varphi_{\epsilon} = 0$ for $t \geq T_0$ and passing to limit $\epsilon \to 0$, we deduce $V(t) f = 0$ for $t \geq T_0.$ This completes the proof of Theorem 4.2. \hfill\qed\\

\begin{corollary}
\label{cor:noeigen}
 With the assumptions of Theorem $4.2$, the operator $G_b$ has no outgoing eigenfunctions in $\Hp$ and $G_b^*$ has no incoming eigenfunctions
in $\Hp$.
\end{corollary}

{\it Proof.}
If $G_b f = \lambda f$ with $f \perp D_{-}^a$ and $u(t, x) = V(t)f = e^{\lambda t}f,\: \Re \lambda < 0$, Theorem 4.2 says that $u(t, x)$ must be disappearing and this yields $f = 0.$ Notice that we can apply only a part of the argument of the proof of Theorem 4.2 leading to $f = 0$ for $|x| \geq b \geq 2 \rho.$ Since $Q f =0,$ we have $(G^{4d} + Q^2- \lambda^{4d})f = 0$. The symbol $ A(\xi)^{4d} +  Q^2(\xi)$ of $G^{4d} + Q^2$ is elliptic with constant coefficients
so  $ f = 0$ in $\Omega.$

For the operator $G_b^*$ a similar argument shows that $G_b^*$ has no incoming eigenfunctions in $H_b^{\perp}.$ 
\hfill\qed

\vskip.2cm

To examine the completeness of the wave operators $W_{\pm}$, notice that  (see \cite{G0})
$$
\overline{{\rm Ran}\: W_{\pm}} 
\ =\
 \Hp \ominus {\mathcal H}_{\pm},
 $$
where
$$
{\mathcal H}_{+} = \{f \in H:\: \lim_{t \to  +\infty} V(t)f = 0\},
\quad
{\mathcal H}_{-} = \{f \in H:\: \lim_{t \to +\infty} V^*(t)f = 0\}.
$$

Next prove  that ${\mathcal H}_{\pm} \perp D_{\pm}^a,\: a \geq \rho.$ In fact, 
if $g \in D_{-}^a$, then $V(t) U_0(-t)g = g$ and
for $f \in {\mathcal H}_{-}$,
$$
(f, g) 
\ =\
 (f, V(t)U_0(-t)g) 
 \ =\
  (V^*(t)f, U_0(-t)g) 
  \ \rightarrow\
   0
   \quad{\rm as}\ \ 
   t \to + \infty.
  $$
A similar argument works for ${\mathcal H}_{+}$. 

\begin{corollary} If $G_b$ has at least one eigenvalue in $\Re z < 0$, then the wave operators $W_{\pm}$ are not complete.
\end{corollary}
{\it Proof.} If $W_{\pm}$ are compete, we have $ {\mathcal H}_{-} = {\mathcal H}_{+}$. The existence of an eigenfunction $f \neq 0$ of $G_b$ with  eigenvalue in $\{\Re \lambda < 0\}$ yields $f \in {\mathcal H}_{-}$ and by the above argument we conclude that $f$ must be disappearing, hence $f = 0$.  
This improves the result in \cite{P3}.

\begin{example} The above argument shows that if $f$ is an eigenfunction for which $G_b f = \lambda f$ with $\Re \lambda < 0$ the outgoing component of $f$ vanishes, so $f$ is always incoming. Our explicit construction 
of such an eigenvalue 
for the Maxwell system in \cite{CPR} with strictly dissipative boundary condition
\begin{equation} \label{eq:4.2}
E_{tan} (1 + \epsilon) - \nu \wedge B_{tan} = 0
\quad
 {\rm on}
 \quad
  |x| = 1,\ \  \epsilon > 0
\end{equation} 
illustrates this situation.
\end{example}

\begin{remark} The space $\Hp \ominus D_{-}^a$ is invariant with respect to $V(t).$ An application of Theorem 4.2 proves that the semigroup
$V(t)$ restricted to this space has no eigenvalues in $\Re z < 0.$ It is an open problem to find dissipative boundary conditions such that
$V(t)$ has no eigenvalues in $\Hp.$
\end{remark}

Thanks to Corollary \ref{cor:noeigen}, we 
may apply the results of Section 5 in \cite{LP2}. For the reader's convenience
 we use the same notations as in \cite{LP2}. Let $\rho_0(B)$ be the component of the resolvent set of $B$ containing right half plane. Then Theorem 4.2 in \cite{LP2} says that $\Sc(z)$ can be  continued analytically from the lower half plane into $-\ii \rho_0(B)$. Thus $\Sc(z)$ is meromorphic in $\C$. The poles of the scattering matrix in $\Im z > 0$ are called scattering resonances and they form a discrete set. The crucial point of our argument is the following.
 
 \begin{definition}
 $\hat{z}$ is a {\bf zero of $\Sc(z)$} if there exists $h \not = 0$ such that $S(\hat{z})h = 0.$
\end{definition} 
 
\begin{theorem}[Theorem 5.6 in \cite{LP2}] 
If $G_b|_{\Hp}$  has no outgoing 
eigenfunctions and $G_b^*|_{\Hp}$ 
has no incoming eigenfunctions, then
 the point spectrum of $G_b\big|_{\Hp}$ is 
 of the form $\ii z$, where $\Im z > 0$ and $z$ is a zero of $\Sc(z)$ or possibly a resonance.
\end{theorem}

\vskip.2cm

{\it Proof of Theorem $\ref{thm:main}$.} For $\Re z < 0$ we know that the index of $G_b - z$ is finite. To prove that this index is 0, it suffices to find at least one point $\hat{z},\: \Re \hat{z} < 0,$ such that $G_b - \hat{z}$ is invertible. 
Equation \eqref{eq:3.4} and the argument of Section 3 show that $\|K^{\#}\|<1$  on a neighborhood of 0, so 
 $\Sc(z)$ has no zeros in that neighborhood.
  Since $B$ has no spectrum on $\ii \R$, there are no resonances on $\R$. By Theorem 4.8  there exists a neighborhood ${\mathcal U}$ of 0, such that for every $z \in {\mathcal U} = \{ z \in \C:\: |z| < \epsilon_0,\: \Im z > 0\}$ the point $\ii z$ is not an eigenvalue of $G_b$.

 An analogous  argument 
 works for $G_b^*$. First  define a scattering operator $S_1 = W_1 \circ W_{+}$ related to $V^*(t)$ and $U_0(t)$. 
 Then introduce a scattering matrix $\Sc_1(z)$ and obtain an analogue of Theorem 4.2. Conclude that there exists a small neighborhood ${\mathcal U}_1$ of 0 such that for $z \in {\mathcal U}_1 = \{z \in \C: \: |z| < \epsilon_1, \: \Im z > 0\}$ the point $\ii z$ is not an eigenvalue of $G_b^*$. Consequently, for $z \in {\mathcal U} \cap {\mathcal U}_1$ the index of $G_b - \ii z$ is 0 and we can apply Theorem 2.4. This proves that outside a discrete set in $\Re z < 0$ the operator $G_b - z$ is invertible. If $z, \: \Re z < 0,$ is not an eigenvalue of $G_b$, but $z \in \sigma(G_b)$, then $z$ must be in the residual spectrum of $G_b$. Then $\dim \Ker (G_b - z) = 0$ and ${\rm codim} \: {\rm Ran}\: (G_b - z) > 0.$ This leads to index $(G_b - z) > 0$ and we obtain a contradiction. The spectrum of $G_b$ in $\Re z < 0$ is formed only by isolated eigenvalues with finite multiplicities.

Finally, the space $\Dp$ is invariant with respect to the semigroup $V(t)$. Thus the generator $G_b$ is the extension of the generator $G_{+}$ of the semigroup $V(t) \vert_{\Dp} = U_0(t)\vert_{\Dp}$. By using the translation representation $\Rc_n$, it is easy to see that $G_{+}$ has spectrum in $\{z \in \C: \Re z \leq 0\}$, so by a well known result the boundary $\ii \R$ of ${\Re z < 0}$ must be included in $\sigma(G_b)$. Since $G_b$ over $\Hp$ has no eigenvalues on $\ii \R$, we deduce that $\ii \R$ is included in the continuous spectrum of $G_b.$ This completes the proof of Theorem 1.4. 
\hfill\qed

Our argument implies easily also the following
\begin{proposition} The operator $G_b$ has no a sequence of eigenvalues $\{z_j\}$ with $\Re z_j < 0, \forall j \in \N$ and $\lim_{j \to \infty} z_j = \ii z_0 \in \ii \R.$
\end{proposition}
{\bf Proof.} Assume that $\{z_j\}_{j \in \N}$ form a sequence of eigenvalues of $G_b$ with $\Re z_j < 0$ and let $\lim z_j = \ii z_0,\: z_0 > 0$ (The case $z_0 < 0$ is completely similar). We will show that in a small neighborhood of $[0, \ii z_0]$ in $\C$ one has only a finite number of eigenvalues. Consider the
interval $J_0 = ]-\epsilon, z_0 + \epsilon[ \subset \R$ with a sufficiently small $\epsilon > 0.$ The scattering matrix ${\mathcal S}(z)$ has no poles on $J_0$, so ${\mathcal S}(z)$ is analytic in a small open neighborhood $W_0 \subset \C$ of $J_0$ and ${\mathcal S}(z) = I + {\mathcal K}(z)$ is  
analytic function in $W_0$ with values Hilbert-Schmidt operators. Since ${\mathcal S}(0)$ is invertible, the analytic Fredholm theorem implies that ${\mathcal S}(z)$ is invertible in $W_0$ outside a discrete set which could have accumulation points only  on the boundary $\partial W_0$. Thus in another neighborhood $W_1 \subset W_0$ of $J_0$ we have at most a finite number points where ${\mathcal S}(z)$ is not invertible. Combining this with Theorem 4.10, one concludes that in $\ii W_1$ we may have at most a finite number of eigenvalues $z$ with $\Re z < 0.$
\hfill\qed

\section{Perturbations}

In \cite{CPR} we constructed an example of an eigenvalue of $G_b$
in $\Hp$ for Maxwell's equations on the exterior of a ball and with 
strictly dissipative boundary conditions.   The construction relied in an
essential way on spherical symmetry.   One could imagine that these
eigenvalues are very sensitive and would disappear under small perturbations of the 
boundary or the boundary conditions.   In this section we show that in fact the
eigenvalues are stable thus extending the construction to the non symmetric
case.

Restrict  attention to the 
 Maxwells equations
   in $\Omega \subset \R^3$. Let $\Lambda(x)$ be a smooth $(2 \times 6)$ matrix-valued function defined for $x \in \partial \Omega$ such that ${\rm rank}\: \Lambda(x) = 2,\: \Lambda(x) u = 0 \:\Leftrightarrow\: u \in \cn (x).$ 
 Assume that the boundary condition $\Lambda(x) = 0$ is strictly maximal dissipative,
 that is,
 $\langle A(\nu(x)) u(x), u(x) \rangle = 0$ for $x \in \partial \Omega$
  implies $u(x) \in \Ker A(\nu(x)).$ 
  
  It is easy to see that $(\Ker G_b)^{\perp} \subset \{ u \in H^1(\Omega):\: \dive u = 0\}.$ For maximally strictly dissipative boundary conditions and $u \in H^1(\Omega)$ we have the following coercive estimate stronger than $(H)$
\begin{equation} \label{eq:5.1}
\|u\|_{H^1(\Omega)}
\ \le\ 
C\Bigl(\| (G - z) u\|_{L^2(\Omega)}
\ +\ 
\|\dive u\|_{L^2(\Omega)} 
\ +\ 
\|\Lambda u\|_{H^{1/2}(\partial\Omega)}
\ +\ 
\|u\|_{L^2(\Omega)}\Bigr)
\,
\end{equation}
with a constant $C > 0$ depending on $\Omega,\: \Lambda$ and $z$.
This estimate is a consequence of the
fact that we can associate to our problem a boundary problem for a second order elliptic system which satisfies the Lopatinski condition (see, for instance, Theorem 2 and Section 5 in \cite{M}).

The perturbation argument needs the following
characterization of eigenvalues.  When
\eqref{eq:5.3}  is  satisfied, it is also satisfied
after  small perturbations of the operator and the boundary condition.
Combined with changes of variables
close to the identity it also applies to small
perturbations of the domain $\Omega$. 

\begin{theorem}  For a fixed $z \in \C,\:\Re z < 0,$
and strictly dissipative boundary condition $\Lambda(x) u = 0$
the following conditions are equivalent:
 
{\bf (i)}  $G_b - z$ is invertible.

{\bf (ii)}  There is a constant $c > 0$ so that for all
$u\in D(G_b) \cap (\Ker G_b)^{\perp}$ we have
\begin{equation}
\label{eq:5.2}
c\|u\|_{H^1(\Omega)}\ \le\  \| (G_b - z)u\|_{L^2(\Omega)}\,.
\end{equation}

{\bf (iii)}  There is a constant $c>0$ so that for
all $u\in H^1(\Omega)$ we have
 \begin{equation}
 \label{eq:5.3}
c \|u\|_{H^1(\Omega)}
 \ \le\
 \| (G - z)u\|_{L^2(\Omega)}
  \ +\ 
  \|\dive u\|_{L^2(\Omega)}  \ +\ \|\Lambda u\|_{H^{1/2}(\partial\Omega) }.
\end{equation}
\end{theorem}

{\it Proof.}   {\bf (i)}$\Leftrightarrow${\bf (ii)}.  Since 
 $G_b - z$ is Fredholm of index zero,
{\bf (i)} 
holds if and only   if there is a constant $c > 0$ so that 
$$
c\|u\|_{D(G_b)} \ \le \ \|(G_b -z)u\|_{L^2(\Omega)}\,,
$$
where $\|u\|_{D(G_b)} = \Bigl( \|u\|^2_{L^2(\Omega)} + \|G u\|^2 _{L^2(\Omega)}\Bigr)^{1/2}$ is the graph norm.

The equivalence of {\bf (i)} and {\bf (ii)}
then follows from the fact that the
graph norm and the $H^1(\Omega)$ norm
are equivalent on $D(G_b) \cap (\Ker G_b)^{\perp}$  and that in turn follows from the condition $(H).$

${\bf (iii)}\:\Rightarrow\:{\bf (ii)}$.
This 
is immediate since for $ u \in D(G_b) \cap (\Ker G_b)^{\perp}$ we have $\Lambda u=0$ and $\dive u=0$.  

  ${\bf (i)}\:\Rightarrow\: {\bf (iii)}$.  
  Assume that {\bf (iii)} is violated. Then since
   \eqref{eq:5.1} is satisfied,
  one can choose a sequence
  $\{u^n\}$ bounded in $H^1(\Omega)$
  with 
  $$
  \|u^n\|_{L^2(\Omega)}=1
  \quad
  {\rm and}
  \quad
  \|(G - z)u^n\|_{L^2(\Omega)}
+
\|\dive u^n\|_{L^2(\Omega)} 
 +
\|\Lambda u^n\|_{H^{1/2}(\partial\Omega)}
\ \to\
0\,.
$$
Passing to a subsequence, which we denote again by $\{u^n\}$, we may suppose that
$\{u^n\}$ converges weakly in $H^1(\Omega)$ to  a limit $w \in H^1(\Omega)$.

{\it The key argument is to show that $\{u^n\}$ is
precompact in $L^2(\Omega)$.
}  Assuming that, one concludes that
there exists a subsequence $\{u^{n_k}\}$  
that converges strongly 
in $L^2(\Omega)$ so $\|w\|_{L^2}=1$.

The weak limits imply that
$ (G - z) w = 0,\:\dive w =0$, and
$\Lambda w=0$ so $w\ne 0$ is an element of 
the kernel of $G_b - z$.
This contradicts ${\bf (i)}$.

It remains to prove the precompactness of 
$\{u^n\}$ in $L^2(\Omega)$. Choose $\phi\in C^\infty(\R^3)$ so that
$\phi=0$ on a neighborhood of $\overline\Omega$
and $\phi=1$ for $|x|\ge R$ for some $0<R<\infty$.
Since $\{u^n\}$ is bounded in $H^1(\Omega)$
and $(1-\phi)$ has compact support, it follows that
$(1-\phi)u_n$ is precompact in $L^2(\R^3)$.
It remains to prove precompactness of 
$\phi u^n$.
Define $u^n:=(E^n,B^n)$ and
$$
zE^n -\curl B^n := f^n_1,
\quad
zB^n +\curl E^n := f^n_2,
\quad
\dive E^n := g^n_1,
\quad
\dive B^n := g^n_2\,.
$$
Setting $f^n = (f_1^n, f_2^n),\: g^n = (g_1^n, g_2^n),$ by construction we have
$$
\lim_{n\to \infty}\ \ 
\| f^n\|_{L^2(\Omega)} = \lim_{n \to \infty} \|g_n\|_{L^2(\Omega)}\ \ =\ 0\,.
$$
Taking the divergence of the equation  $(z - G)\Bigl(\begin{matrix} E_n\\ B_n\end{matrix}\Bigl) = \Bigl(\begin{matrix} f^n_1\\f^2_n \end{matrix}\Bigr)$
shows that
$$
zg^n_1= z\,\dive E^n = \dive f^n_1\,,
\quad
 zg^n_2= z\,\dive B^n = \dive f^n_2\,,
 $$
 so 
 $$
\lim_{n\to \infty}\ \ 
\|\,  \dive f^n\, \|_{L^2(\Omega)}\ \ =\ 0\,.
$$
Furthermore,
\begin{equation}
\label{eq:5.4}
\begin{aligned}
z\phi E^n -\curl \phi B^n = \phi f^n_1 -[\curl,\phi]B^n&:= \widetilde f^n_1\,,
\cr
z\phi B^n +\curl \phi E^n = \phi f^n_2 + [\curl,\phi] E^n&:=\widetilde f^n_2\,,
\cr
\dive \phi E^n \ =\ \phi g^n_1 + [\dive,\phi]E^n&:= \widetilde g^n_1\,,
\cr
\dive \phi B^n \ =\ \phi g^n_2 + [\dive,\phi]B^n&:=\widetilde g^n_2\,.
\end{aligned}
\end{equation}
Since $f^n$ and $g^n$ tend to zero
in $L^2(\Omega)$, and $E^n$ and $B^n$ are
bounded in $H^1(\Omega)$, it follows from the expressions
defining them that $\widetilde f^n$ and $\widetilde g^n $ 
are precompact in $L^2(\R^3)$.

We prove the precompactness of $\{\phi E^n\}$.  The case of 
$\phi B^n$ is entirely analogous.  Compute
\begin{equation} \label{eq:5.5}
\begin{aligned}
z^2 \phi E^n&= \curl z\phi B^n + z\widetilde f^n_1 =-\curl\curl(\phi E^n) + \curl \widetilde f^n_2+ z\widetilde f^n_1
\cr
&=
\Delta (\phi E^n) + \grad\dive (\phi E^n)+\curl \widetilde f^n_2
+ z\widetilde f^n_1.
\end{aligned}
\end{equation}
Therefore
$$
(z^2-\Delta)(\phi E^n)
\ =\ 
\grad \widetilde g^n_1+\curl \widetilde f^n_2
+ z\widetilde f^n_1\,.
$$
Inverting $z^2-\Delta$ on tempered distributions on $\R^3$ by
using the fundamental solution, yields
$$
\phi E^n
 \ =\ 
 \frac{e^{z|x| } }{4\pi |x|}
\,*\,
\Big(
\grad \widetilde g^n_1+\curl \widetilde f^n_2
+ z\widetilde f^n_1
\Big)\,.
$$
The expression in parentheses is precompact
in $H^{-1}(\R^3)$.  It follows that
$\phi E^n$ is precompact in $H^1(\R^3)$
and therefore in $L^2(\Omega)$. This completes
the proof.
\hfill\qed

\begin{remark} We can generalize Theorem 5.1 for systems with strictly dissipative boundary conditions satisfying the condition $(E)$ for which we have an analogue of (\ref{eq:5.1}).  Then  $\dive u = 0$ is replaced by $Qu = 0$ and for the proof we may exploit the 
fundamental solution of $z^2 - G^2 + Q^* Q.$
\end{remark}

It follows from the characterization (5.3)
that if  $G_b - z$ is invertible, then
the same is
true after small perturbations of $G$, $\Lambda$, $z$ and 
the domain $\Omega$.
 If $|z-\hat z|=r$ lies in $\Re z<0$ and
there are no eigenvalues on this circle, the same is true
after small perturbation, and the spectral projector 
$$
\frac1{2\pi i}\,
\oint_{|z-\hat z|=r}   (z-G_b)^{-1}\ dz
$$
also depends continuously on perturbations.
 After perturbation, an eigenvalue
 may split into a finite number of points
(no larger than the finite rank of the associated spectral
projector).  



\end{document}